\documentclass[12pt, a4paper]{article}

\usepackage[latin1]{inputenc}
\usepackage{german}
\usepackage{multicol}

\usepackage{a4wide} 

\begin{document}
\begin{center}
  \large\textbf{Der rechnende Dichter}\normalsize\\[2ex]
(Vortrag f\"ur die Studienstiftung des deutschen Volkes am 16. Februar
2001 in M\"unchen)\\[2ex]

\textbf{Christian Blohmann}\\[1ex]
International University Bremen, School of Engineering and Science\\
Campus Ring 1, D-28759 Bremen\\

\end{center}


    
Wenn ich heute Abend im Haus eines Geisteswissenschaftlers die
Erlaubnis habe, über Mathematik zu reden, zwar über sozusagen
literarisch domestizierte, aber eben doch über Mathematik, so bin ich
mir der damit verbundenen Verantwortung wohl bewusst. Nur wenig Themen
scheinen besser geeignet, die Stimmung einer Abendgesellschaft auf den
Nullpunkt zu bringen, als mathematische.

So hatte Goethe in seinem Haus wohl ein diesbezügliches Gespächsverbot
erlassen, sagte er doch von einem seiner Gäste: "`Er bringt das
allerfremdeste, was in mein Haus kommen kann, die Mathematik an meinen
Tisch; wobey wir jedoch schon eine Convention geschlossen haben, daß
nur im alleräußersten Falle von Zahlen die Rede seyn
darf."'\footnote{Johann Wolfgang Goethe, Werke, Weimarer Ausgabe, Bd.
  IV 20, Weimar 1887-1919, S. 224.}

Auch Herr Broich stellte sich im Vorbereitungsgespräch zu diesem Abend
besorgt hinter Goethes Geselligkeitspostulat und wies mich darauf hin,
auch ja darauf zu achten, dass mein Vortrag "`verständlich"', d.h.
nicht zu mathematisch sein dürfe. Diesen Sorgen möchte ich gleich
vorab entgegentreten. Wenn ich heute Abend über den "`rechnenden
Dichter"' reden will, so muss man weder einen
mathematisch-physikalischen, noch einen
literatur\emph{wissenschaftlichen} Vortrag befürchten. Mir, als
Angehörigem der, wie Georg Lukács sagt, desanthropomorphisierenden,
also entmenschlichenden exakten Wissenschaften,
leistet seit dem Beginn meines Studiums die schöngeistige Lektüre
insofern einen großen Dienst, als sie mir hilft, meiner
Entmenschlichung etwas entgegen zu setzen. Leider liefert die
Literatur Naturwissenschaftlern und Mathematikern nur selten
tatsächliche Identifikationsangebote, ja, der Drang der Literaten,
sich mit den exakten Wissenschaften zu beschäftigen, scheint nur
gering ausgeprägt zu sein, womöglich, weil, wie uns Musil wissen läßt,
die Seele eines rechten Schöngeistes das ist, "`was sich verkriecht,
wenn man von algebraischen Reihen hört."'\footnote{Robert Musil, Der
  Mann ohne Eigenschaften, hrsg. v. Adolf Frisé, Reinbek bei Hamburg
  1978, S. 103.}

Wird die Grenzüberschreitung dennoch gewagt, sei es durch Verwendung
mathematischer Begriffe als Metaphern, durch Übertragung
mathematischer oder physikalischer Methoden in die Poetologie, durch
Einsatz mathematisch gebildeter Romanhelden, durch literarische
Reflexionen auf empirische Naturerkenntnis oder nur, um den Leser
durch vermeintliches Geheimwissen einzuschüchtern, meistens ist
das Ergebnis bemerkenswert. Zum Beleg möchte ich eine kleine,
unsystematische und auch \emph{nicht} chronologisch angeordnete, persönliche
Auswahl meiner liebsten Fundstellen von der Klassik über die Moderne
bis hin zur Postmoderne anführen und ein wenig kommentieren.

Zurück zu Goethe. Er war der Mathematik nicht generell abgeneigt und
hatte selbst Privatunterricht bei einem Jenaer Mathematikprofessor
genommen. Darüber schreibt er an Frau von Stein: "`Algebra ist
angefangen worden, sie macht noch ein grimmig Gesicht [...]"'. Zwei
Tage später: "`Vielleicht treffe ich irgendwo eine Lücke durch die ich
mich einschleiche."' Doch nach weiteren zwei Tagen: "`[...] ich werde
es zu meinem Wesen nicht brauchen können, da das Handwerck ganz ausser
meiner Sphäre liegt"'\footnote{Goethe, a.a.O., Bd. IV 7, S. 219ff.}
Liegt die Mathematik aber innerhalb der eigenen Sphäre, so wünscht man
sich als Naturwissenschaftler oder Mathematiker oft nichts sehnlicher
als eine Gelegenheit, der Begeisterung etwa über einen gefundenen
Beweis Ausdruck zu verleihen. Um dieses zugegebenermaßen oft
unangebrachte Mitteilungsbedürfnis von Mathematikern in Gesellschaft
nachvollziehen zu können, scheint wie bei Heinrich von Kleist das eine
oder andere mathematische Erfolgserlebnis vonnöten.  Dieser schreibt:
\begin{quote}
  "`Aber wenn ich einen mathematischen Lehrsatz ergründet habe, dessen
  Erhabenheit und Größe mir auch die Seele füllte, wenn ich nun mit
  diesem Eindruck in eine Gesellschaft trete, wem darf ich mich
  mitteilen, wer versteht mich?"'\footnote{Heinrich von Kleist,
    Sämtliche Werke und Briefe, hrsg. v. Helmut Sembdner, Bd. 2,
    München 1984, S. 497f.}
\end{quote}
Haben umgekehrt persönliche Misserfolge das Verhältnis zur Mathematik
bestimmt, ist oft eine generelle Ablehnung die Folge.
Alfred Döblin geht soweit zu behaupten:
\begin{quote}
  "`Eine lächerliche Sache überhaupt, diese Mathematik auf den
  Schulen. Für die meisten wertlos, ein abseitiges Gedankenspiel, eine
  Qual, weil ohne Anschauung, ohne Ziel, ohne Bindung mit einem Leben.
  Man soll diese Art Abstraktion verbieten oder in die Akademien
  schicken."'\footnote{Alfred Döblin, Schriften zu Leben
    und Werk, hrsg. v. Erich Kleinschmidt, Freiburg 1986, S. 152.}
\end{quote}
Dagegen scheint Heinrich Heine immerhin mit einer der Grundrechenarten
seinen Frieden gemacht zu haben:
\begin{quote}
  "`Ich aber hatte in der Schule meine Not mit den vielen Zahlen! Mit
  dem eigentlichen Rechnen ging es noch schlechter. Am besten begriff
  ich das Subtrahieren, und da gibt es eine sehr praktische
  Hauptregel: \glq Vier von drei geht nicht, da muß ich eins
  borgen\grq -- ich rate aber jedem, in solchen Fällen immer einige
  Groschen mehr zu borgen; denn man kann nicht
  wissen."'\footnote{Heinrich Heine, Werke in 4 Bänden, Frankfurt am
  Main 1994, Bd. 2, S. 197.}
\end{quote}
Doch was macht die Mathematik mit den Menschen, wenn sie sich ihr
nicht verweigern, sondern, und sei es auch nur zeitweilig,
verschreiben? Inwiefern kann man, über die Schwierigkeiten, Mathematik
zu kommunizieren, hinausgehend, tatsächlich von einer Entmenschlichung,
einer Beeinträchtigung des menschlichen Lebens durch Mathematik
sprechen? Mathematische Beschäftigung wird immer wieder als leidvoll
beschrieben, als Tä\-tig\-keit, die dem Gefühlsleben nicht gerecht wird
bzw. dieses schmerzhaft beeinträchtigt. So litt etwa Heinrich von
Kleist trotz oder wegen seiner mathematischen Erfolge:
\begin{quote}
  "`Bei dem ewigen Beweisen und Folgern verlernt das Herz fast zu
  fühlen; und doch wohnt das Glück nur im Herzen, nur im Gefühl, nicht
  im Kopfe, nicht im Verstande. Das Glück kann nicht, wie ein
  mathematischer Lehrsatz bewiesen werden, es muß empfunden werden,
  wenn es dasein soll."'\footnote{Kleist, a.a.O., S. 494.}
\end{quote}
Wie muss es erst denjenigen ergehen, die sich ihr Leben lang mit
Mathematik beschäftigen? Da sie selbst in der Regel nicht oder nur
schlecht dichten -- abschreckende Beispiele gibt es durchaus --, ist
man auf literarische Schilderungen von Nichtmathematikern angewiesen.
Hans Magnus Enzensberger, der sich zeit seines Lebens mit Mathematik
und Mathematikern beschäftigt hat, verdanken wir ein teilweise
erschreckend zutreffendes Gedicht über \emph{Die Mathematiker}, in dem
er mit dem alltagssprachlichen Doppelsinn mathematischer Fachtermini
spielt:
\begin{verse}
"`Wurzeln, die nirgends wurzeln,\\
Abbildungen für geschlossene Augen,\\
Keime, Büschel, Faltungen, Fasern:\\
diese weißeste aller Welten\\
mit ihren Garben, Schnitten und Hüllen\\
ist euer Gelobtes Land.

Hochmütig verliert ihr euch\\
im Überabzählbaren, in Mengen\\
von leeren, mageren, fremden\\
in sich dichten und Jenseits-Mengen.

Geisterhafte Gespräche\\
unter Junggesellen:\\
die Fermatsche Vermutung,\\
der Zermelosche Einwand,\\
das Zornsche Lemma.

Von kalten Erleuchtungen\\
schon als Kinder geblendet,\\
habt ihr euch abgewandt,\\
achselzuckend,\\
von unseren blutigen Freuden.

Wortarm stolpert ihr,\\
selbstvergessen,\\ 
getrieben vom Engel der Abstraktion,\\
über Galois-Felder und Riemann-Flächen,\\
knietief im Cantor-Staub,\\
durch Hausdorffsche Räume.

Dann, mit vierzig, sitzt ihr,\\
o Theologen ohne Jehova,\\
haarlos und höhenkrank\\
in verwitterten Anzügen\\
vor dem leeren Schreibtisch,\\
ausgebrannt, o Fibonacci,\\
o Kummer, o Gödel, o Mandelbrot,\\
im Fegefeuer der Rekursion."'\footnote{Hans Magnus Enzensberger, Zukunftsmusik,
  Frankfurt am Main 1993, S. 26f.} 
\end{verse}


Vorherrschende Abneigung, grundlegendes Misstrauen oder Desinteresse
gegenüber der Mathematikerzunft halten die Schriftsteller aber nicht
davon ab, sich mathematisch-natur\-wissen\-schaft\-li\-cher Metaphern,
Vergleiche oder Bilder zu bedienen. Bei Goethe ist es nicht nur die
Physik, sondern auch die Mathematik, die er, trotz der erwähnten
Vorbehalte, mit hörbarem Respekt herbeizitiert:
\begin{verse}
"`Nicht meinem Witz ward solche Gunst beschert,\\
Zwei Götterschwestern haben mich belehrt:\\
Physik voran, die jedes Element\\
Verbinden lehrt, wie sie es erst getrennt;\\
Das Unwägbare hat für sie Gewicht,\\
Und aus dem Wasser lockt sie Flammenlicht,\\
Läßt Unbegreifliches dann sichtbar sein\\
Durch Zauberei im Sondern, im Verein.

Doch erst zur Tat erregt den tiefsten Sinn\\
Geometrie, die Allbeherrscherin:\\
Sie schaut das All durch ein Gesetz belebt,\\
Sie mißt den Raum und was im Raume schwebt;\\
Sie regelt streng die Kreise der Natur,\\
Hiernach die Pulse deiner Taschenuhr;\\
Sie öffnet geistig grenzenlosen Kreis\\
Der Menschenhände kümmerlichstem Fleiß."'\footnote{Goethe, Schriften
  zur Naturwissenschaft, Leopoldina-Ausgabe, Weimar 1947ff., Bd. I 2,
  S. 368f.} 
\end{verse}
Hier wird die streng logisch-empirische Gültigkeit der exakten
Wissenschaften in allegorischer Weise zum Sinnbild einer über den
Menschen, ja sogar über die Götter hinausweisenden Ordnungsmacht, der
"`Allbeherrscherin"'. Wegen dieses tranzendentalen Potenzials
erfreuen sich einige mathematische Metaphern großer Beliebtheit, wie
die Vorstellung von zwei parallelen Geraden, die sich nirgends, oder,
wenn man so sagen will, im Unendlichen schneiden. Etwa in Günter Eichs
Gedicht \emph{Des Toten gedenken}, in dem das Nirgendwo, an dem sich zwei
Parallelen schneiden, eine Jenseitsvorstellung repräsentiert:
\begin{verse}
"`Dorthin gehen,\\
wo die Parallelen sich schneiden.\\
Die Forderungen der Logik\\
durch Träume erfüllen."'\footnote{Günter Eich, Gesammelte Werke,
  hrsg. v. Axel Vieregg, Frankfurt am Main 1991, Bd. 1, S. 154.} 
\end{verse}
Oder in humoristischer Form in Morgensterns \emph{Die zwei
  Parallelen}, in dem die Parallelen den Weg in eine paradiesische
  Unendlichkeit markieren.
\begin{quote}
\begin{multicols}{2}
\setcounter{unbalance}{3}
\noindent
"`Es gingen zwei Parallelen\\
ins Endlose hinaus,\\
zwei kerzengerade Seelen\\
und aus solidem Haus.

\noindent
Sie wollten sich nicht schneiden\\
bis an ihr seliges Grab:\\
das war nun einmal der beiden\\
geheimer Stolz und Stab.

\noindent
Doch als sie zehn Lichtjahre\\  
gewandert neben sich hin,\\
da wards dem einsamen Paare\\
nicht irdisch mehr zu Sinn.

\noindent
War'n sie noch Parallelen?\\
Sie wußten's selber nicht. --\\
sie flossen nur wie zwei Seelen\\
zusammen durch ewiges Licht.

\noindent
Das ewige Licht durchdrang sie,\\
da wurden sie eins in ihm;\\
die Ewigkeit verschlang sie,\\
als wie zwei Seraphim."'\footnote{Christian Morgenstern, Alle
  Galgenlieder, Wiesbaden 1947, S. 247.}
\end{multicols}
\end{quote}
Doch aus bereits genannten Gründen sind es oft negative Assoziationen,
die durch mathematische Metaphorik evoziert werden sollen. So wird
Heine in den \emph{Reisebildern} durch die Begegnung mit einer offenbar
unangenehmen Person an das ihm nicht minder unangenehm erscheinende
Repertoire schulmathematischer Begrifflichkeit erinnert:
\begin{quote}
  "`Doktor Saul Ascher mit seinen abstrakten Beinen, mit seinem engen,
  transzendentalgrauen Leibrock, und mit seinem schroffen, frierend
  kalten Gesichte, das einem Lehrbuche der Geometrie als Kupfertafel
  dienen konnte. Dieser Mann [...] war eine personifizierte gerade
  Linie."'\footnote{Heine, Werke, a.a.O., Bd. 2, S. 11.}
\end{quote}
Um eine tiefsinnigere Nutzung des mathematischen Metaphernsteinbruchs
bemüht zeigt sich Wedekind, der von Mathematik nicht viel verstand,
der aber mathematische Konzepte sehr gekonnt im literarischen
Zusammenhang zu instrumentieren wusste.  In der Kindertragödie
\emph{Frühlings Erwachen} ist es nicht zuletzt die Mathematik, an der
die pubertierenden, von den Erwachsenen unverstandenen Protagonisten
zu leiden und teilweise auch zu sterben haben: Auf Marthas Nachfrage,
weshalb sich Moritz umgebracht habe, sagt Ilse "'Parallelepipedon!"'.
Gemeint ist Moritz' Scheitern am entsprechenden mathematischen
Schulstoff. Einzig Melchior gelingt es, den Übertritt ins
Erwachsenenleben mehr oder minder unbeschadet zu überstehen.
Unterstützung erfährt er dabei von dem "`vermummten Herrn"', einer
rätselhaften Figur, dessen Rolle Wedekind, sofern es ihm möglich war, selbst
gespielt hat, und die es versteht, mathematische Metaphorik gekonnt
für eigene Zwecke, in diesem Fall die Veranschaulichung
moralphilosophischer Zusammenhänge, zu handhaben.
\begin{quote}
  "`Melchior: Wie denken Sie über Moral? [...] \\
  Der vermummte Herr: Unter Moral verstehe ich das reelle Produkt
  zweier imaginärer Größen. Die imaginären Größen sind \emph{Sollen}
  und \emph{Wollen}. Das Produkt heißt Moral und läßt sich in seiner
  Realität nicht leugnen."'\footnote{Frank Wedekind, Frühlings
    Erwachen, Stuttgart 1971, S. 69.}
\end{quote}
Das Setzen der Rechnungseinheit $i = \sqrt{-1}$ als zunächst
abstraktes gedankliches Konstrukt, das im Raum der herkömmlichen
Zahlen nicht existiert und in diesem Sinne eingebildet, also imaginär
ist, aber doch zum Verständnis der "`realen"' Zahlen hilfreich ist,
wird verglichen mit der moralphilosophischen Erklärung des
Zustandekommens einer ebenso realen sittlichen Handlung durch erstens
das Wissen um moralische Prinzipien -- Sollen -- und zweitens die
Bereitschaft diesen zu folgen -- Wollen.

In seinem Gedicht \emph{Von der Algebra der Gefühle} geht Enzensberger
noch einen Schritt weiter, indem er die menschlichen Gefühle in
Hilberts Hotel wohnen läßt. David Hilbert, womöglich der
einflussreichste Mathematiker der Neuzeit, hatte zur Veranschaulichung
eines mengentheoretischen Paradoxons folgendes Gleichnis herangezogen:
In einem besonders großen Hotel, in dem es zu jeder natürlichen Zahl
-- 1,2,3,... usw. -- genau ein Zimmer mit dieser Zimmernummer gibt,
sind alle Zimmer belegt. Ein neuer Gast kommt an. Der Portier bittet
kurzerhand den Gast aus Zimmer 1 in Zimmer 2 umzuziehen, den aus
Zimmer 2 in Zimmer 3, aus 3 in 4 und so fort, bis ins Aschgraue, wie
Enzensberger sagen würde. In das freigewordene Zimmer 1 kann nun der
Neuankömmling einziehen. Wiederum sind alle Zimmer belegt. Da sich die
Anzahl der Zimmer nicht geändert hat, kann sich bei voller Belegung
auch die Anzahl der Gäste nicht geändert haben, obwohl ein neuer Gast
dazugekommen ist... das Leben ist seltsam in Hilberts Hotel. Hier nun
das Gedicht:
\begin{verse}
"`Ich habe oft das Gefühl (brennend,\\
dunkel, undefinierbar usw.),\\
daß das Ich keine Tatsache ist,\\
sondern ein Gefühl,\\
das ich nicht loswerde.\pagebreak

Ich hege es, lasse ihm freien Lauf,\\
erwidere es, von Fall zu Fall.\\
Aber es ist nur eins unter vielen.

Die Menge der Gefühle ist abzählbar unendlich,\\
d.h. sie lassen sich im Prinzip numerieren,\\
bis ins Aschgraue.

Die Nummer der Eifersucht\\
ist offensichtlich die Sieben.\\
Auch die Angst ist prim.\\
Und ich habe das dumpfe Gefühl,\\
daß die Demütigung die 188 auf ihrer Stirn trägt --\\
eine Zahl ohne Eigenschaften.

Auch das Gefühl, numeriert zu sein,\\
ist vermutlich längst numeriert,\\
nur wozu und von wem?

Das erhabne Gefühl des Zorns\\
bewohnt ein anderes Zimmer\\
in Hilberts Hotel\\
als das Gefühl,\\
über den Zorn erhaben zu sein.

Und nur wer sich hingeben kann\\
dem abstrakten Gefühl\\
für die Abstraktion, der weiß,\\
daß es in manchen sehr hellen Nächten\\
den Wert $i$ anzunehmen pflegt.

Dann wieder läuft es mir kalt\\
über den Rücken, das Gefühl,\\
ein Paket zu sein,\\
das gefühllose, pelzige Gefühl,\\
von dem die Zunge zu bersten droht\\
nach der Injektion,\\
wenn sie dem Zahn auf den Zahn fühlt,\\
oder die Peinlichkeit\\
mit ihrem durchdringenden Bleigeschmack,\\
das mächtige Gefühl der Ohnmacht,\\
das unaufhaltsam der Null zustrebt,\\
und das falsche Gefühl\\
der wahren Empfindung\\
mit seinem abscheulichen Kettenbruch.

Dann erfüllt mich\\
eine Schnittmenge aus gemischten Gefühlen,\\
schuldig, fremd, wohl, verloren,\\
alles auf einmal.

Nur dem höchsten der Gefühle\\
wäre das Ich nicht gewachsen.\\
Statt Aufwallungen zu suchen\\
mit dem Limes $\infty$,\\
läßt es sich lieber\\
eine Minute lang übermannen\\
vom Schauder des eisig heißen Wassers\\
unter der Dusche, dessen Nummer\\
noch keiner entziffert hat."'\footnote{Enzensberger, Kiosk. Neue
  Gedichte, Frankfurt am Main 1995, S. 47ff.}
\enlargethispage{1ex}
\end{verse}
Heinrich Heine hätte sich wohl gegen diesen Versuch, die Zahlen in den
menschlichen Gefühlen und damit in den Menschen selbst zu sehen, verwehrt:
\begin{quote}
  "`Da war ein trübseliger Minister, respektabler Bankier, guter
  Hausvater, guter Christ, guter Rechner ... und vor lauter Zahlen sah
  er weder die Menschen noch ihre drohenden Mienen ... Wahrlich, es
  ist töricht, wenn man nur die Personen sieht in den Dingen, so ist
  es noch törichter, wenn man in den Dingen nur die Zahlen sieht. Es
  gibt aber Kleingeister, die aufs pfiffigste beide Irrtümer zu
  verschmelzen suchen, die sogar in den Personen die Zahlen suchen,
  womit sie uns die Dinge erklären wollen."'\footnote{Heine, Werke,
  a.a.O., Bd. 3, S.196f.}
\end{quote}


Bereits von antiken Mathematikern wurde diejenige Methode entwickelt,
die heute als die mathematische Methode schlechthin gilt, nämlich das
streng logische Schließen ausgehend von einigen wenigen postulierten
Grundannahmen, den Axiomen, auf komplexe mathematische Sachverhalte.
Diese erstmals in Euklids \emph{Elementen} überlieferte Technik des
deduktiven Beweises, war lange Zeit in einen Dornröschenschlaf
versunken, bis sie in der Neuzeit von Descartes wiederentdeckt wurde,
der in seinen \emph{Regeln zur Ausrichtung der Erkenntniskraft}
fordert, "`daß, wer den richtigen Weg zur Wahrheit sucht, mit keinem
Gegenstand umgehen darf, über den er nicht eine den arithmetischen
oder geometrischen Beweisen gleiche Gewißheit gewinnen
kann."'\footnote{René Descartes, Regeln zur Ausrichtung der
  Erkenntniskraft, hrsg. v. Lüder Gäbe, Hamburg 1979, S.~8.} Als
erstes von Spinoza in seiner \emph{Ethica more geometrico demonstrata}
(Ethik nach geometrischer Methode dargestellt) aufgegriffen,
verbreitete sich die deduktive Methode schnell in allen Bereichen der
Wissenschaft. In der 1758 erschienenen, äußerst wirkungsmächtigen
\emph{Histoire des Mathématiques} -- auch in Goethes Bibliothek befand
sich ein Exemplar -- kolportiert Jean Étienne Montucla eine Anekdote,
derzufolge der Mathematiker Roberval nach der Lektüre von Racines
Iphigenie provokant gefragt haben soll: "`Qu'est-ce que cela prouve?"'
-- Was beweist das schon?

Den hingeworfenen Fehdehandschuh ergreift in Deutschland als einer der
ersten
Johann Carl Wezel. In seinem Roman \emph{Lebensgeschichte Tobias
  Knauts} ironisiert er:
\begin{quote}
  "`Meinen kürzesten, deutlichsten, bündigsten Beweis will ich, wie in
  allen Sachen, also auch hier gebrauchen. In gehöriger Form steht er
  also:\\
  1. Mir ist es unbegreiflich\\
  2. - - - \\
  3. - - - \\
  4. - - - und so ins unendliche fort.\\
  Den will ich doch sehen, der wider diesen Beweis etwas einzuwenden
  weis! Vielleicht könnten einige schwergläubige Leute, denen die
  Wahrheit niemals Wahrheit ist, wenn sie nicht in dem nemlichen
  Kleide erscheint, in welchem sie sie alle Tage sehen -- Vielleicht
  könnten solche, sage ich, bey der vorhergehenden Deducktion
  gravitätisch sich beym Kinne fassen, und mit einer verachtenden
  vielbedeutenden Mine denken oder sagen: Qu'est-ce que cela
  prouve?"'\footnote{Johann Karl Wezel, Lebensgeschichte Tobias
    Knauts, des Weisen, sonst der Stammler genannt, Bd. 1, Stuttgart
    1971, S.18f.}
\end{quote}
Andererseits gibt es den ein wenig verstiegen Versuch Friedrich
Schlegels, für den "`Mathematik = Magie"'\footnote{Friedrich Schlegel,
  Kritische Ausgabe, hrsg. v. Ernst Behler, München u.a. 1958ff., Bd.
  19, S. 10.}  gilt, auf der Basis mathematischer Begrifflichkeit eine
Poetologie zu entwerfen. Mit den Formelhaften Abkürzungen $F$, $S$ und
$M$ für Fantastik, Sentimentalität und Mimik sieht das dann so aus:
\begin{quote}
  "`$\frac{F}{0}$ $\frac{S}{0}$ $\frac{M}{0}$ sind
  die \emph{poetischen Ideen}.\\
  Das poetische Ideal $=\sqrt[1/0]{\frac{FSM^{1/0}}{0}}=$
  Gott"'\footnote{Ebd., Bd. 16, S. 148.}\\
\end{quote}
Gegen einen derartigen mathematischen Einschüchterungsversuch, mag man
sich -- der Anachronismus sei mir gestattet -- mit Lessing wehren:
\begin{quote}
  "`Hu! würde ein stolzer Algebraist murmeln, Ihr mein Freund ein
  Philosoph? Laßt einmal sehen. Ihr versteht doch wohl einen
  hyperbolischen Afterkegel zu kubieren? Oder nein -- Könntet Ihr eine
  Exponential-Größe differentieren? Es ist eine Kleinigkeit: hernach
  wollen wir unsre Kräfte in was größern versuchen. Ihr schüttelt den
  Kopf? Nicht? Nu da haben wirs. Bald wollte ich wetten, Ihr wißt
  nicht einmal, was eine Irrational-Größe ist? Und werft Euch zu einem
  Philosoph auf? O Verwegenheit! o Zeit! o
  Barbarei!"'\footnote{Gotthold Ephraim Lessing, Werke, hrsg. v.
    Herbert Göpfert, München 1970-79, Bd. 3, S. 689f.}
\end{quote}
Während Schlegel den Poeten gerne zum Mathematiker erzöge, deklariert
Novalis, der selbst über eine fundierte mathematische Ausbildung
verfügte, den Mathematiker kurzerhand zum Poeten:    
\begin{quote}
  "`[...] wenn der Mathematiker wircklich etwas richtiges thut, so
  thut ers, als poetischer philosoph. [...] Der poetische Philosoph
  ist en état de Createur absolu. Ein Kreis, ein Triangel werden schon
  auf diese Art creirt. Es kommt ihnen nichts zu, als was der
  Verfertiger ihnen zukommen läßt etc."'\footnote{Novalis, Schriften,
    hrsg. v. Richard Samuel, Darmstadt 1977-88, Bd. 3, S. 415.}
\end{quote}
Novalis sieht in der gemeinsamen Schöpfungsfreiheit die Verwandtschaft
von Poesie und Mathematik begründet. \emph{Hier} ist es die zunächst freie
Erschaffung eines poetische Kunstwerks, das nur den intrinsischen
Regeln ästhetischer Stimmigkeit genügen muss, \emph{dort} die mathematische
Setzung von Axiomen und Definitionen, die lediglich logisch
wiederspruchsfrei zu sein hat. Aus heutiger wissenschaftstheoretischer
Sicht ist diese Verwandtschaft durchaus begründet, ist es doch der
empirische Wahrheitsbegriff, der die Naturwissenschaften und einige
Sozialwissenschaften von den Geisteswissenschaften aber eben auch von der
Mathematik trennt.

Franz Grillparzer sieht im Deutschland seiner Zeit durch eine solche
Verbrüderung von Poesie und Mathematik zu einem Geist deutscher
Gründlichkeit die Dichtung in ihrer ästhetischen Wirksamkeit bedroht.
Der gelehrte deutsche Theaterbesucher überfrachtete nach Grillparzer
-- angespornt durch die Literaturkritik -- die Dichtung mit
literaturfremden Erwartungen und verbaute sich damit auf seine
gründliche Art das eigentliche literarische Erlebnis. In seinem 1850
erschienenen Gedicht \emph{Gründlichkeit} heißt es:
\pagebreak

\begin{multicols}{2}
\setcounter{unbalance}{3}
\noindent
"`Wie viel, im Reich des Geistes gar,\\
hängt ab von Ort und Zeit\\
Was falsch einst, gilt uns heut für wahr,\\
Für dumm, was sonst gescheit.
\\[1ex]
\noindent
Und mancher, den die eigne Zeit\\
Verspottet und verlacht,\\
Lebt' er in unsern Tagen, heut,\\
Sein Glück wär' längst gemacht.
\\[1ex]
\noindent
So jener Mathematikus\\
Im heiteren Paris,\\
Setzt ins Theater nie den Fuß,\\
Da Zahlen nur gewiß.
\\[1ex]
\noindent
Doch einst die Freunde brachten ihn\\
Ins Schauspielhaus mit Glück,\\
Man gab ein Schauspiel von Racine,\\
Des Meisters Meisterstück.
\\[1ex]
\noindent
Da wird denn rings Begeistrung laut,\\
Man weint, man klatscht, man tobt,\\
Was man gehört, was man geschaut,\\
Wird \emph{eines} Munds gelobt.
\\[1ex]
\noindent
Nur unser Mathematikus\\
Sah stieren Augs das Spiel,\\
Als ihn der Freunde Schar am Schluß\\
Befragt: wie's ihm gefiel,

\noindent
Ob ihn ergriff der Dichtung Macht,\\
Des Unglücks Jammerruf?\\
Doch er erwidert mit Bedacht:\\
\glq Mais qu'est-ce que cela prouve?\grq
\\[1ex]
\noindent
Da tönt Gelächter rings umher,\\
Das Wort durchläuft die Stadt,\\
Und ein Jahrhundert oder mehr\\
Lacht sich die Welt nicht satt.
\\[1ex]
\noindent
O armer Mann, du kamst zu früh,\\
Und nicht am rechten Ort;\\
In unsers Deutschlands Angst und Müh'\\
Erkennt man erst dein Wort.
\\[1ex]
\noindent
Wo man Ideen nur begehrt,\\
Von Glut und Reiz entfernt,\\
Man, bis zum Halse schon gelehrt,\\
Noch im Theater lernt --
\\[1ex]
\noindent
Dort ruft ein jeder Kritikus,\\
Was auch der Dichter schuf,\\
Wie jener Mathematikus:\\
\glq Mais qu'est-ce que cela prouve?\grq "'\footnote{Franz Grillparzer,
  Gedichte, hrsg. v. August Sauer, Wien 1932, S. 222f.}
\end{multicols}


Mathematik und in ihrem Gefolge die Naturwissenschaften sind im
Verlauf des 20. Jahrhunderts zu Leitwissenschaften geworden, die
unsere Welt nachhaltig beeinflußt und verändert haben, mit
Ausstrahlungskraft in nahezu jeden gesellschaftlichen Bereich. Einst
gefeiert als die treibende Kraft des Fortschritts geriet das
mathematische Denken mit der um sich greifenden Fortschrittsskepsis
jedoch zunehmend unter kulturkritischen und kulturpessimistischen
Druck. 

Die -- um mit Max Weber zu sprechen -- "`Entzauberung der Welt"' durch
technologische Veränderungen hatte seine Schattenseiten, die in der
Fortschrittseuphorie zunächst übersehen wurden, die aber, waren diese
Schattenseiten erst einmal offen zu Tage getreten, den exakten
Wissenschaften angelastet wurden, man könnte sagen: in der
Erleichterung darüber, einen Sündenbock gefunden zu haben.  In der
Regel resultieren diese, auch heute noch verbreiteten Vorwürfe aus
illusionären Ansprüchen, die an die Wissenschaften herangetragen
werden und die geradezu notwendig zu Enttäuschungen führen.

Die Literatur des 20. Jahrhunderts ist voll von Versuchen, die
Mathematik in ihre Schranken zu weisen, ihr die Deutungsmacht für
nicht-exakte, menschliche Lebensbereiche zu entziehen. So etwa bei
Volker Braun in dem Gedicht \emph{Mitteilung an die reifere Jugend}
von 1965:
\pagebreak

\begin{verse}
  "`Nein, die Bäume unserer Lust könnt ihr nicht konstruieren.\\
  Zirkel und Lineal geben ein Plakat ab, nicht das Leben.\\
  Unsere Vergnügen sind zu toll und gefährlich, um sie in Serien zu\\
  produzieren.\\
  Unser Glück ist total: es läßt sich nicht ausrechnen.\\
  Unsere Vorsicht vor uns ist vergeblich: wir sind
  maßlos."'\footnote{Volker Braun, Provokation für mich. Gedichte, Halle
    1965, S. 19}
\end{verse}

Doch was kann Mathematik für die Moderne leisten, wenn sie, befreit
von falschen Heilserwartungen, auf das Innovationspotenzial ihrer
Methode hin befragt wird? Ich möch\-te mit einem längeren Zitat enden,
und zwar einer Passage aus Musils Mann ohne Eigenschaften, nicht ohne
Grund die heimliche Bibel mancher literarisch interessierter
Naturwissenschaftler. Musil nimmt Bezug auf die Urträume der
Menschheit, die, technologisch realisiert, ihrer zauberischen
Anziehungskraft beraubt erscheinen, verharrt allerdings nicht bei
dieser kulturpessimistischen Position, sondern verweist auf den im
mathematischen Denken sich manifestierenden
"`Möglichkeitssinn"'\footnote{Musil, Mann ohne Eigenschaften, a.a.O.,
  S. 16.} als die Fähigkeit, sich hypothetisch und schöpferisch auf
die Wirklichkeit zu beziehen:

\begin{quote}
  "`[...]  es ist nicht zu leugnen, daß [...] diese Urtäume nach
  Meinung der Nichtmathematiker mit einemmal in einer ganz anderen
  Weise verwirklicht waren, als man sich das ursprünglich vorgestellt
  hatte. Münchhausens Posthorn war schöner als die fabrikmäßige
  Stimmkonserve, der Siebenmeilenstiefel schöner als ein Kraftwagen,
  [...] die Vögel zu verstehn, schöner als eine tierpsychologische
  Studie über die Ausdrucksbewegungen der Vogelstimme. Man hat
  Wirklichkeit gewonnen und Traum verloren. [...]  Man braucht
  wirklich nicht viel darüber zu reden, es ist den meisten Menschen
  heutzutage ohnehin klar, daß die Mathematik wie ein Dämon in alle
  Anwendungen unseres Lebens gefahren ist. [...] Und so hat es auch
  schon damals, als Ulrich Mathematiker wurde, Leute gegeben, die den
  Zusammenbruch der europäischen Kultur voraussagten, weil kein
  Glaube, keine Liebe, keine Einfalt, keine Güte mehr im Menschen
  wohne, und bezeichnender Weise sind sie alle in ihrer Jugend
  schlechte Mathematiker gewesen. Damit war später für sie bewiesen,
  daß die Mathematik, Mutter der exakten Naturwissenschaft, Großmutter
  der Technik, auch Erzmutter jenes Geistes ist, aus dem schließlich
  Giftgase und Kampfflieger aufgestiegen sind. [...]  Von Ulrich [dem
  Protagonisten des Romans] dagegen konnte man mit Sicherheit das eine
  sagen, daß er die Mathematik liebte, wegen der Menschen, die sie
  nicht ausstehen mochten. [...] Er sah, daß sie in allen Fragen, wo
  sie sich für zuständig hält, anders denkt als gewöhnliche Menschen.
  [...] Es geht in der Wissenschaft so stark und unbekümmert und
  herrlich zu wie in einem Märchen. Und Ulrich fühlte: die Menschen
  wissen das bloß nicht; sie haben keine Ahnung, wie man schon denken
  kann; wenn man sie neu denken lehren könnte, würden sie auch anders
  leben."'\footnote{Musil, Der Mann ohne Eigenschaften, a.a.O., S.
    39f.}
\end{quote}

\end{document}